\documentclass[12pt,twoside]{article}

\usepackage{amsmath}
\usepackage{amssymb}
\usepackage{amsthm}

\date{}

\begin{document}

\centerline{}

\centerline{}

\centerline {\Large{\bf Fractional Optimal Control in the Sense of
Caputo}}

\centerline{}

\centerline{\Large{\bf and the Fractional Noether's
Theorem\footnote{Int. Math. Forum (accepted 11 Dec 2007)}}}

\centerline{}

\centerline{\bf {Gast\~{a}o S.~F.~Frederico}}

\centerline{}

\centerline{Department of Science and Technology}

\centerline{University of Cape Verde}

\centerline{Praia, Santiago, Cape Verde}

\centerline{\texttt{silves36@hotmail.com}}

\centerline{}

\centerline{\bf {Delfim F.~M.~Torres}}

\centerline{}

\centerline{Department of Mathematics}

\centerline{University of Aveiro}

\centerline{3810-193 Aveiro, Portugal}

\centerline{\texttt{delfim@ua.pt}}

\newtheorem{theorem}{\quad Theorem}[section]

\newtheorem{definition}[theorem]{\quad Definition}

\newtheorem{corollary}[theorem]{\quad Corollary}

\newtheorem{lemma}[theorem]{\quad Lemma}

\newtheorem{example}[theorem]{\quad Example}

\newtheorem{remark}[theorem]{\quad Remark}

\begin{abstract}
The study of fractional variational problems with derivatives in
the sense of Caputo is a recent subject, the main results being
Agrawal's necessary optimality conditions of Euler-Lagrange and
respective transversality conditions. Using Agrawal's
Euler-Lagrange equation and the Lagrange multiplier technique, we
obtain here a Noether-like theorem for fractional optimal control
problems in the sense of Caputo.
\end{abstract}




\section{Introduction}

Noether's theorem, published in 1918 \cite{Noether:1918}, is a
central result of the calculus of variations that explains all
physical laws based upon the action principle. It is a very
general result, asserting that ``to every variational symmetry of
the problem there corresponds a conservation law''. Noether's
principle gives powerful insights from the various transformations
that make a system invariant. For instance, in mechanics the
invariance of a physical system with respect to spatial
translation gives conservation of linear momentum; invariance with
respect to rotation gives conservation of angular momentum; and
invariance with respect to time translation gives conservation of
energy \cite{NoetherLivroFr,CD:delfimCPAA04}. The calculus of
variations is now part of a more vast discipline, called optimal
control \cite{Pontryagin}, and Noether's principle still holds in
this more general setting \cite{CD:delfimEJC,Torres06}.

Fractional derivatives play an increasing role in mathematics,
physics and engineering \cite{Agrawal:2004b,CD:Hilfer:2000,%
Kilbas,CD:MilRos:1993,CD:Podlubny:1999,CD:SaKiMa:1993}
and the theory of the calculus of variations has been extended in
order to deal with more general systems containing non-integer
derivatives \cite{CD:Agrawal:2002,CD:Agrawal:2004a,CD:Agrawal:2006,%
AvkarBaleanu,Cresson2007,El-NabulsiMMAS,Klimek:2005,Zbl05083320,CD:Riewe:1996}.
The new fractional variational calculus provide a more realistic
approach to physics \cite{Baleanu05,Baleanu06,Cresson05,%
Klimek:2002,CD:Muslih:2005,CD:Riewe:1997}, permitting to consider
nonconservative systems in a natural way---a very important issue
since closed systems do not exist: forces that do not store
energy, so-called nonconservative or dissipative forces, are
always present in real systems. Nonconservative forces remove
energy from the systems and, as a consequence, the standard
conservation laws cease to be valid. However, it is still possible
to obtain a Noether-type theorem which covers both conservative
and nonconservative cases \cite{CD:Djukic:1980,comGastaoParis06}.
In a more general way, formulations of Noether's theorem were
proved for fractional problems with left and right
Riemann-Liouville fractional derivatives
\cite{Cresson2007,CD:FredericoTorres:2007,CD:FredericoTorres:2006,NonlinearDyn}.

In \cite{CD:Agrawal:2006} Agrawal proves a version of the
Euler-Lagrange equations for fractional problems of the calculus
of variations in the sense of Caputo. One of the interesting
aspects of the new theory is that both Caputo and
Riemann-Liouville derivatives play a role in Agrawal's
Euler-Lagrange equations. Here we use the results of
\cite{CD:Agrawal:2006} to formulate a Noether-type theorem in the
general context of the fractional optimal control in the sense of
Caputo.


\section{Fractional Derivatives}
\label{sec:fdRL-C}

In this section we collect the well-known definitions of
fractional derivatives in the sense of Riemann-Liouville and
Caputo (\textrm{see}
\cite{CD:Agrawal:2002,CD:Agrawal:2006,CD:MilRos:1993,CD:Podlubny:1999,CD:SaKiMa:1993}).

\begin{definition}[Fractional derivative in the sense of Riemann-Liouville]
Let $f$ be an integrable continuous function in the interval
$[a,b]$. For $t \in [a,b]$, the left Riemann-Liouville fractional
derivative $_aD_t^\alpha f(t)$ and the right Riemann-Liouville
fractional derivative $_tD_b^\alpha f(t)$, of order $\alpha$, are
defined by
\begin{gather}
_aD_t^\alpha f(t) =
\frac{1}{\Gamma(n-\alpha)}\left(\frac{d}{dt}\right)^{n}
\int_a^t (t-\theta)^{n-\alpha-1}f(\theta)d\theta \, , \label{eq:DFRLE} \\
_tD_b^\alpha f(t)
=\frac{1}{\Gamma(n-\alpha)}\left(-\frac{d}{dt}\right)^{n}
\int_t^b(\theta - t)^{n-\alpha-1}f(\theta)d\theta \, ,
 \label{eq:DFRLD}
\end{gather}
where $n \in \mathbb{N}$, $n-1 \leq \alpha < n$, and $\Gamma$ is
the Euler gamma function.
\end{definition}

\begin{definition}[Fractional derivative in the sense of Caputo]
Let $f$ be an integrable continuous function in $[a,b]$. For $t
\in [a,b]$, the left Caputo fractional derivative $_a^CD_t^\alpha
f(t)$ and the right Caputo fractional derivative $_t^CD_b^\alpha
f(t)$, of order $\alpha$, are defined in the following way:
\begin{gather}
_a^CD_t^\alpha f(t) = \frac{1}{\Gamma(n-\alpha)} \int_a^t
(t-\theta)^{n-\alpha-1}\left(\frac{d}{d\theta}\right)^{n}
f(\theta)d\theta \, , \label{eq:DFCE} \\
_t^CD_b^\alpha f(t) =\frac{1}{\Gamma(n-\alpha)} \int_t^b(\theta -
t)^{n-\alpha-1}\left(-\frac{d}{d\theta}\right)^{n}f(\theta)d\theta
\, , \label{eq:DFCD}
\end{gather}
where $n \in \mathbb{N}$, $n-1 \leq \alpha < n$.
\end{definition}

\begin{remark}
If $\alpha\in \mathbb{\mathbb{N}}$, equalities
\eqref{eq:DFRLE}-\eqref{eq:DFCD} give the classical derivatives
\begin{equation*}
\begin{split}
 _aD_t^\alpha f(t) &=\, _a^CD_t^\alpha f(t)
 =\left(\frac{d}{dt}\right)^\alpha f(t) \, , \\
 _tD_b^\alpha f(t) &=\,_t^CD_b^\alpha f(t)=
\left(-\frac{d}{dt}\right)^\alpha f(t) \, .
\end{split}
\end{equation*}
\end{remark}

\begin{remark}
The Caputo fractional derivative of a constant is always equal to
zero. This is not the case with the Riemann-Liouville fractional
derivative.
\end{remark}


\section{Main Results}

Our main result is a Noether-type theorem for fractional optimal
control problems in the sense of Caputo
(Theorem~\ref{theo:tndfc}). As a corollary, we obtain a Noether
theorem for the fractional problems of the calculus of variations
(Corollary~\ref{theo:tnfc2}).

The fractional optimal control problem in the sense of Caputo is
introduced, without loss of generality, in Lagrange form:
\begin{gather}
\label{Pcap} I[q(\cdot),u(\cdot)] = \int_a^b
L\left(t,q(t),u(t)\right) dt
\longrightarrow \min \, , \tag{$P_{C}$}\\
_a^CD_t^\alpha q(t)=\varphi\left(t,q(t),u(t)\right) \, , \notag
\end{gather}
where functions $q : [a,b] \rightarrow \mathbb{R}^n$ satisfy
appropriate boundary conditions. The Lagrangian $L :[a,b] \times
\mathbb{R}^{n}\times \mathbb{R}^{m} \rightarrow \mathbb{R}$ and
the velocity vector $\varphi:[a,b] \times \mathbb{R}^{n}\times
\mathbb{R}^m\rightarrow \mathbb{R}^{n}$ are assumed to be
functions of class $C^{1}$ with respect to all their arguments. We
also assume, without loss of generality, that $0<\alpha\leq1$. In
conformity with the calculus of variations, we are considering
that the control functions $u(\cdot)$ take values on an open set
of $\mathbb{R}^m$. Along the work we denote by $\partial_{i}L$, $i
= 1,2,3$, the partial derivative of function
$L(\cdot,\cdot,\cdot)$ with respect to its $i$th argument.

\begin{definition}[Process]
An admissible pair $(q(\cdot),u(\cdot))$ which satisfies the
control system $_a^CD_t^\alpha
q(t)=\varphi\left(t,q(t),u(t)\right)$ of problem \eqref{Pcap} is
said to be a \emph{process}.
\end{definition}

\begin{remark}
Choosing $\alpha=1$, Problema~\eqref{Pcap} is reduced to the
classical problem of optimal control theory \cite{Pontryagin}:
\begin{gather}
\label{eq:coc}
 I[q(\cdot),u(\cdot)] = \int_a^b L\left(t,q(t),u(t)\right) dt
\longrightarrow \min \, , \\
\dot{q}(t)=\varphi\left(t,q(t),u(t)\right) \, . \notag
\end{gather}
\end{remark}

\begin{remark}
The fundamental fractional problem of the calculus of variations
in the sense of Caputo, first introduced in \cite{CD:Agrawal:2006}
\begin{equation}
\label{eq:pfccv1} I[q(\cdot)] = \int_a^b
L\left(t,q(t),\,_a^CD_t^\alpha q(t)\right) \longrightarrow \min \, ,
\end{equation}
is a particular case of \eqref{Pcap}: we just need to choose
$\varphi(t,q,u)=u$.
\end{remark}

The fractional Hamiltonian formalism introduced in
\cite{El-NabulsiMMAS} is easily adapted to our present context.
Using the standard Lagrange multiplier technique, we rewrite
problem \eqref{Pcap} in the following equivalent form:
 \begin{equation}
\label{eq:Pcap1} I[q(\cdot),u(\cdot),p(\cdot)] = \int_a^b
\left[{\mathcal
H}\left(t,q(t),u(t),p(t)\right)-p(t)\cdot\,_a^CD_t^\alpha
q(t)\right]dt \longrightarrow \min \, ,
\end{equation}
where the Hamiltonian ${\mathcal H}$ is defined by
\begin{equation}
\label{eq:Hnd} {\mathcal H}\left(t,q,u,p\right) =L\left(t,q,u\right)+p
\cdot \varphi\left(t,q,u\right) \, .
\end{equation}

\begin{remark}
In the context of classical mechanics, $p$ is interpreted as the
generalized momentum. In the optimal control literature, the
multiplier $p$ is known as the adjoint variable.
\end{remark}

We now proceed with the usual steps for obtaining necessary
optimality conditions in the calculus of variations. We begin by
computing the variation $\delta I$ of functional \eqref{eq:Pcap1}:
\begin{equation}
\label{eq:PcapV} \delta I=\int_a^b\left[\partial_2{\mathcal
H}\cdot\delta q+\partial_3{\mathcal H}\cdot\delta u+\partial_4{\mathcal
H}\cdot\delta p-\delta p\cdot \,_a^CD_t^\alpha q-p\cdot
\delta\left(_a^CD_t^\alpha q\right)\right]dt\, ,
\end{equation}
where $\delta q$, $\delta u$ and $\delta p$ are the variations of
$q$, $u$, and $p$ respectively. Using the results obtained in
\cite{CD:Agrawal:2002,CD:Agrawal:2004a,CD:Agrawal:2006}, equation
\eqref{eq:PcapV} is equivalent to
\begin{multline*}
 \delta I=\int_a^b\left[\left(\partial_2{\mathcal
H}-\,_tD_b^\alpha p\right)\cdot\delta q+\partial_3{\mathcal
H}\cdot\delta u+\left(\partial_4{\mathcal H}- \,_a^CD_t^\alpha
q\right)\cdot\delta p\right]dt \\
-\left.\left(_tD_b^{\alpha-1} p\right)\cdot\delta q\right|_a^b\,
\end{multline*}
where the fractional derivatives of $p(t)$ are in the sense of
Riemann-Liouville (in contrast with the fractional derivative of
$q(t)$ which is taken in the sense of Caputo). Standard arguments
conduce us to the following result.

\begin{theorem}
\label{theo:pmpdfc} If $(q(\cdot),u(\cdot))$ is an optimal process
for problem \eqref{Pcap}, then there exists a function
$p(\cdot)\in C^{1}([a,b];\mathbb{R}^{n})$ such that for all $t\in
[a,b]$ the tuple $(q(\cdot),u(\cdot),p(\cdot))$ satisfy the
following conditions:
\begin{itemize}
\item the Hamiltonian system
\begin{equation*}
\begin{cases}
\partial_{2} {\mathcal H}(t, q(t), u(t),p(t))=\,_tD_b^{\alpha} p(t) \, , \\
 \partial_{4} {\mathcal H}(t, q(t), u(t),
p(t))=\,_a^CD_t^\alpha q(t) \, ;
\end{cases}
\end{equation*}
\item the stationary condition
\begin{equation*}
 \partial_{3} {\mathcal H}(t, q(t), u(t), p(t))=0 \, ;
\end{equation*}
\item the transversality condition
\begin{equation*}
\left(_tD_b^{\alpha-1} p\right)\cdot\delta q|_a^b=0\, ;
\end{equation*}
\end{itemize}
with ${\mathcal H}$ given by \eqref{eq:Hnd}.
\end{theorem}

\begin{definition}
\label{fracPE}
A triple $(q(\cdot),u(\cdot),p(\cdot))$ satisfying
Theorem~\ref{theo:pmpdfc} will be called \emph{a fractional
Pontryagin extremal}.
\end{definition}

\begin{remark}
For the fundamental fractional problem of the calculus of
variations in the sense of Caputo \eqref{eq:pfccv1} we have ${\mathcal
H}=L+p\cdot u$. It follows from Theorem~\ref{theo:pmpdfc} that
\begin{equation}
\label{eq:dedELfPMP}
\begin{gathered}
\,{_a^CD_t^\alpha} q=u\, ,\\
_tD_b^\alpha p=\partial_2 L\, ,\\
\partial_3 {\mathcal H}=0 \Leftrightarrow p=-\partial_3 L\Rightarrow
{_tD_b}^\alpha p=-_tD_b^\alpha \partial_3 L\,.
\end{gathered}
\end{equation}
Comparing both expressions for $_tD_b^\alpha p$, we arrive to the
fractional Euler-Lagrange equations proved by O. P. Agrawal
\cite[\S 3]{CD:Agrawal:2006}:
\begin{equation}
\label{eq:E-Lca}
\partial_2 L+\,_tD_b^\alpha \partial_3 L=0\, .
\end{equation}
In other words, for the problem of the calculus of variations
\eqref{eq:pfccv1} the fractional Pontryagin
extremals give Agrawal's Euler-Lagrange extremals.
\end{remark}

\begin{remark}
Our optimal control problem \eqref{Pcap} only involves Caputo
fractional derivatives but both Caputo and Riemann-Liouville
fractional derivatives appear in the necessary optimality
condition given by Theorem~\ref{theo:pmpdfc}. This is different
from \cite{CD:Agrawal:2004a,CD:FredericoTorres:2006,NonlinearDyn}
where the necessary conditions only involve the same type of
derivatives (Riemann-Liouville) as those in the definition of the
fractional optimal control problem. This fact also occurs in the
particular case of the calculus of variations as noted in
\cite{CD:Agrawal:2006}: the Riemann-Liouville derivative is
present in the necessary condition of optimality \eqref{eq:E-Lca}
but not in the formulation of the problem \eqref{eq:pfccv1}.
\end{remark}

The notion of variational invariance for problem \eqref{Pcap} is
defined with the help of the equivalent problem \eqref{eq:Pcap1}.

\begin{definition}[Invariance of \eqref{Pcap} without transformation of time]
\label{def:inv-coc} We say that functional \eqref{Pcap} is
invariant under the one-parameter family of infinitesimal
transformations
\begin{equation*}
\label{eq:tinffc}
\begin{cases}
\bar{q}(t) = q(t) + \varepsilon\xi(t,q,u,p) + o(\varepsilon) \, ,\\
\bar{u}(t) = u(t) + \varepsilon\varsigma(t,q,u,p) + o(\varepsilon) \, ,\\
\bar{p}(t) = p(t) +\varepsilon\varrho(t,q,u,p) + o(\varepsilon) \, ,\\
\end{cases}
\end{equation*}
if and only if
\begin{multline}
\label{eq:inv-cap} \int_{t_{a}}^{t_{b}}\left[{\mathcal
H}\left(t,q(t),u(t),p(t)\right)-p(t)\cdot\,_a^CD_t^\alpha
q(t)\right]dt \\=\int_{t_{a}}^{t_{b}}\left[{\mathcal
H}\left(t,\bar{q}(t),\bar{u}(t),\bar{p}(t)\right)
-\bar{p}(t)\cdot\,_a^CD_t^\alpha\bar{q}(t)\right]dt\,
\end{multline}
for any subinterval $[{t_{a}},{t_{b}}] \subseteq [a,b]$.
\end{definition}

Noether's theorem will be proved following similar steps as those
in \cite{CD:FredericoTorres:2007} for the problems of the calculus
of variations in the Riemann-Liouville sense.

\begin{lemma}[necessary and sufficient condition of invariance]
If functional \eqref{Pcap} is invariant, in the sense of
Definition~\ref{def:inv-coc}, then
\begin{multline}
\label{eq:cnsidfc}
\partial_{2}{\mathcal H}\left(t,q,u,p\right) \cdot \xi
+ \partial_3 {\mathcal H}\left(t,q,u,p\right)\cdot\varsigma+
\left(\partial_{4} {\mathcal H}\left(t,q,u,p\right)-{_a^CD_t^\alpha
q}\right)\cdot\varrho-p\cdot\,_a^CD_t^\alpha \xi = 0 \, .
\end{multline}
\end{lemma}

\begin{proof}
Since condition \eqref{eq:inv-cap} is to be valid for any
subinterval $[{t_{a}},{t_{b}}] \subseteq [a,b]$, we can write
\eqref{eq:inv-cap} in the following equivalent form:
\begin{equation}
\label{eq:inv-cap1} {\mathcal
H}\left(t,q(t),u(t),p(t)\right)-p(t)\cdot\,_a^CD_t^\alpha q(t)
={\mathcal H}\left(t,\bar{q}(t),\bar{u}(t),\bar{p}(t)\right)
-\bar{p}(t)\cdot\,_a^CD_t^\alpha\bar{q}(t) \, .
\end{equation}
We differentiate both sides of \eqref{eq:inv-cap1} with respect to
$\varepsilon$ and then substitute $\varepsilon$ by zero. The
definition and properties of the Caputo fractional derivative
permit us to write that
\begin{multline}
\label{eq:SP1} 0 = \partial_{2} {\mathcal H}\left(t,q,u,p\right)\cdot\xi
+ \partial_{3} {\mathcal H}\left(t,q,u,p\right) \cdot
\varsigma+\left(\partial_{4} {\mathcal
H}\left(t,q,u,p\right)-{_a^CD_t^\alpha
q}\right)\cdot\varrho\\
-p\cdot\frac{d}{d\varepsilon} \left[\frac{1}{\Gamma(n-\alpha)}
\int_a^t
(t-\theta)^{n-\alpha-1}\left(\frac{d}{d\theta}\right)^{n}q(\theta)d\theta\right.
\\+\left.\frac{\varepsilon}{\Gamma(n-\alpha)}\int_a^t
(t-\theta)^{n-\alpha-1}\left(\frac{d}{d\theta}\right)^{n}\xi \,d\theta\right]_{\varepsilon=0}\\
= \partial_{2} {\mathcal H}\left(t,q,u,p\right)\cdot\xi(t,q) +
\partial_{3} {\mathcal H}\left(t,q,u,p\right) \cdot
\varsigma+\left(\partial_{4} {\mathcal
H}\left(t,q,u,p\right)-{_a^CD_t^\alpha q}\right)\cdot\varrho\\
-p\cdot\frac{1}{\Gamma(n-\alpha)}\int_a^t
(t-\theta)^{n-\alpha-1}\left(\frac{d}{d\theta}\right)^{n}\xi\,
d\theta
 \, .
\end{multline}
Expression \eqref{eq:SP1} is equivalent to \eqref{eq:cnsidfc}.
\end{proof}

In this work we propose the new notion of \emph{fractional
conservation law in the sense of Caputo}. For that, we introduce
the operator $\mathcal{D}_{t}^{\omega}$.

\begin{definition}
Given two functions $f$ and $g$ of class $C^{1}$ in the interval
$[a,b]$, we introduce the following operator:
\begin{equation*}
\mathcal{D}_{t}^{\omega}\left[f,g\right] = -g \, {_tD_b^\omega} f +f
\, {_a^CD_t^\omega} g \, ,
\end{equation*}
where $t \in [a,b]$ and $\omega \in \mathbb{R}_0^+$.
\end{definition}

\begin{remark}
Similar operators were used in \cite{El-NabulsiMMAS} and
\cite[Def.~19]{CD:FredericoTorres:2007}. While here
$\mathcal{D}_{t}^{\omega}$ depends both on Caputo and
Riemann-Liouville derivatives, the ones introduced in
\cite{El-NabulsiMMAS,CD:FredericoTorres:2007} involve
Riemann-Liouville derivatives only. We note that in the classical
context $\omega = 1$, and $\mathcal{D}_{t}^{1}[f,g] = g f' + f g'
= \frac{d}{dt}\left(f g\right) = \mathcal{D}_{t}^{1}[g,f]$.
\end{remark}

\begin{definition}[fractional conservation law in the sense %
of Caputo -- \textrm{cf.} Def.~23 of \cite{CD:FredericoTorres:2007}]
\label{def:leicofc}
A quantity $C_f\left(t,q(t),{_a^CD_t^\alpha}q(t),u(t),p(t)\right)$ is said to
be a fractional conservation law in the sense of Caputo if
it is possible to write $C_{f}$ as a sum of products,
\begin{multline}
\label{eq:somaPrdca}
C_{f}\left(t,q(t),{_a^CD_t^\alpha}q(t),u(t),p(t)\right)\\ =
\sum_{i=1}^{r}
C_{i}^1\left(t,q(t),{_a^CD_t^\alpha}q(t),u(t),p(t)\right) \cdot
C_{i}^2\left(t,q(t),{_a^CD_t^\alpha}q(t),u(t),p(t)\right)
\end{multline}
for some $r \in \mathbb{N}$, and for every $i = 1,\ldots,r$ the
pair $C_{i}^1$ and $C_{i}^2$ satisfy one of the following
relations:
\begin{equation}
\label{eq:fracc-cons}
\mathcal{D}_{t}^{\alpha}\left[C_{i}^{1}\left(t,q(t),{_a^CD_t^\alpha}q(t),u(t),p(t)\right)
,C_{i}^{2}\left(t,q(t),{_a^CD_t^\alpha}q(t),u(t),p(t)\right)\right]
= 0
\end{equation}
or
\begin{equation}
\label{eq:fracc1-cons}
\mathcal{D}_{t}^{\alpha}\left[C_{i}^{2}\left(t,q(t),
{_a^CD_t^\alpha}q(t),u(t),p(t)\right),
C_{i}^{1}\left(t,q(t),{_a^CD_t^\alpha}q(t),u(t),p(t)\right)\right]
= 0
\end{equation}
along all the fractional Pontryagin extremals
(Definition~\ref{fracPE}).
\end{definition}

\begin{remark}
If $\alpha=1$ \eqref{eq:fracc-cons} and \eqref{eq:fracc1-cons}
coincide, and $C_{f}$ \eqref{eq:somaPrdca} satisfy the classical
definition of conservation law:
$\frac{d}{dt}\left[C_f(t,q,\dot{q})\right]=0$.
\end{remark}

\begin{lemma}[Noether's theorem without transformation of time]
\label{theo:tnfc1} If functional \eqref{Pcap} is invariant in the
sense of Definition~\ref{def:inv-coc}, then $$p(t)\cdot\xi$$ is a
fractional conservation law in the sense of Caputo.
\end{lemma}

\begin{proof}
We use the conditions of Theorem~\ref{theo:pmpdfc} in the
necessary and sufficient condition of invariance
\eqref{eq:cnsidfc}:
\begin{equation*}
\begin{split}
0=&-\partial_{2}{\mathcal H} \cdot \xi - \partial_3 {\mathcal
H}\cdot\varsigma - \left(\partial_{4} {\mathcal H}-{_a^CD_t^\alpha
q}\right)\cdot\varrho + p\cdot\,_a^CD_t^\alpha \xi\\
=&-\xi\cdot\,_tD_b^\alpha p+p\cdot\,_a^CD_t^\alpha
\xi\\=&\mathcal{D}_t^\alpha\left[p,\xi\right]
 \, .
\end{split}
\end{equation*}
\end{proof}

\begin{definition}[Invariance of \eqref{Pcap}]
\label{def:invnd-co} Functional \eqref{Pcap} is said to be
invariant under the one-parameter infinitesimal transformations
\begin{equation}
\label{eq:tinf}
\begin{cases}
\bar{t} = t + \varepsilon\tau(t,q,u,p) + o(\varepsilon) \, ,\\
\bar{q}(t) = q(t) + \varepsilon\xi(t,q,u,p) + o(\varepsilon) \, ,\\
\bar{u}(t) = u(t) + \varepsilon\varsigma(t,q,u,p) + o(\varepsilon) \, ,\\
\bar{p}(t) = p(t) +\varepsilon\varrho(t,q,u,p) + o(\varepsilon) \, ,\\
\end{cases}
\end{equation}
if and only if
\begin{multline}
\label{eq:invnd-co} \int_{t_{a}}^{t_{b}}\left[{\mathcal
H}\left(t,q(t),u(t),p(t)\right)-p(t)\cdot{_a^CD_t^\alpha q(t)}\right]dt\\
=\int_{\bar{t}(t_a)}^{\bar{t}(t_b)} \left[{\mathcal
H}\left(\bar{t},\bar{q}(\bar{t}),\bar{u}(\bar{t}),\bar{p}(\bar{t})\right)
-\bar{p}(\bar{t})\cdot {_a^CD_{\bar{t}}^{\alpha}
\bar{q}}(\bar{t})\right]d\bar{t}
\end{multline}
for any subinterval $[{t_{a}},{t_{b}}] \subseteq [a,b]$.
\end{definition}

Next result provides an extension of Noether's theorem
\cite{CD:delfimEJC} for fractional optimal control problems in the
sense of Caputo.

\begin{theorem}[Noether's theorem for fractional optimal control problems]
\label{theo:tndfc} If the functional \eqref{Pcap} is invariant
under the one-parameter infinitesimal transformations
\eqref{eq:tinf}, then
\begin{multline}
\label{eq:tn-cap}
C_f\left(t,q(t),{_a^CD_t^\alpha}q(t),u(t),p(t)\right)
\\ =\left[{\mathcal
H}(t,q(t),u(t),p(t))-(1-\alpha)p(t)\cdot{_a^CD_t^\alpha
q(t)}\right\rceil\tau-p(t)\cdot\xi
\end{multline}
is a fractional conservation law in the sense of Caputo (see
Definition~\ref{def:leicofc}).
\end{theorem}

\begin{remark}
If $\alpha=1$ problem \eqref{Pcap} takes the classical form
\eqref{eq:coc} and Theorem~\ref{theo:tndfc} gives the conservation
law of \cite{CD:delfimEJC}:
\begin{equation*}
C(t,q(t),u(t),p(t)) =\left[{\mathcal
H}(t,q(t),u(t),p(t))\right]\tau-p(t)\cdot\xi\, .
\end{equation*}
\end{remark}

\begin{proof}
Every non-autonomous problem \eqref{Pcap} is equivalent to an
autonomous one by artificially considering $t$ as a dependent
variable. For that we consider a Lipschitzian transformation
\begin{equation*}
[a,b]\ni t\longmapsto \sigma f(\lambda) \in [\sigma_{a},\sigma_{b}]
\end{equation*}
satisfying the condition
$t_{\sigma}^{'}=\frac{dt(\sigma)}{d\sigma}=f(\lambda) = 1$ for
$\lambda=0$, such that \eqref{eq:Pcap1} takes the form
\begin{multline*}
\bar{I}[t(\cdot),q(t(\cdot)),u(t(\cdot)),p(t(\cdot))] =
\int_{\sigma_{a}}^{\sigma_{b}} \left[{\mathcal
H}\left(t(\sigma),q(t(\sigma)),u(t(\sigma)),p(t(\sigma)\right)\right.
\\-\left.
p(t(\sigma))\cdot{_{\sigma_{a}}^CD_{t(\sigma)}^{\alpha}q(t(\sigma))}\right]t_{\sigma}^{'}
d\sigma\, ,
\end{multline*}
where $t(\sigma_{a}) = a$, $t(\sigma_{b}) = b$ and
\begin{equation*}
\begin{split}
_{\sigma_{a}}^C&D_{t(\sigma)}^{\alpha}q(t(\sigma))\hspace*{-0.5cm}\\
&= \frac{1}{\Gamma(n-\alpha)} \int_{\frac{a}{f(\lambda)}}^{\sigma
f(\lambda)}\left({\sigma
f(\lambda)}-\theta\right)^{n-\alpha-1}\left(\frac{d}{d\theta}\right)^{n}
q\left(\theta f^{-1}(\lambda)\right)d\theta    \\
&= \frac{(t_{\sigma}^{'})^{-\alpha}}{\Gamma(n-\alpha)}
\int_{\frac{a}{(t_{\sigma}^{'})^{2}}}^{\sigma}
(\sigma-s)^{n-\alpha-1}\left(\frac{d}{ds}\right)^{n}
q(s)ds  \\
&= (t_{\sigma}^{'})^{-\alpha}\,\,{{^C_\chi
D_{\sigma}^{\alpha}q(\sigma)}},\,\left(
\chi={\frac{a}{(t_{\sigma}^{'})^{2}}}\right)\, .
\end{split}
\end{equation*}
Then, we have
\begin{equation*}
\begin{split}
&\bar{I}[t(\cdot),q(t(\cdot)),u(t(\cdot)),p(t(\cdot))] \\
&=\int_{\sigma_{a}}^{\sigma_{b}} \left[{\mathcal
H}\left(t(\sigma),q(t(\sigma)),u(t(\sigma)),p(t(\sigma)\right)
-p(t(\sigma))\cdot(t_{\sigma}^{'})^{-\alpha}\,\,{{^C_\chi
D_{\sigma}^{\alpha}q(\sigma)}}\right]t_{\sigma}^{'}
d\sigma \\
&\doteq \int_{\sigma_{a}}^{\sigma_{b}} \bar{{\mathcal
H}}_{f}\left(t(\sigma),t_{\sigma}^{'},q(t(\sigma)),u(t(\sigma)),p(t(\sigma)
,{{^C_\chi D_{\sigma}^{\alpha}q(\sigma)}}\right)d\sigma \\
& = \int_a^b \left[{\mathcal
H}\left(t,q(t),u(t),p(t)\right)-p(t)\cdot\,_a^CD_t^\alpha
q(t)\right]dt\\
& = I[q(\cdot),u(\cdot),p(\cdot)] \, .
\end{split}
\end{equation*}
If functional $I[q(\cdot),u(\cdot),p(\cdot)]$ is invariant in the
sense of Definition~\ref{def:invnd-co}, then the functional
$\bar{I}[t(\cdot),q(t(\cdot)),u(t(\cdot)),p(t(\cdot))]$ is
invariant in the sense of Definition~\ref{def:inv-coc}. Applying
Lemma~\ref{theo:tnfc1}, we obtain that
\begin{equation}
\label{eq:tempo2}
C_f\left(t(\sigma),t_{\sigma}^{'},q(t(\sigma)),u(t(\sigma)),p(t(\sigma))
,{{^C_\chi
D_{\sigma}^{\alpha}q(\sigma)}}\right)=p(t(\sigma))\cdot\xi+\psi(t(\sigma))\tau
\end{equation}
is a fractional conservation law in the sense of Caputo. For
$\lambda=0$,
\begin{equation}
\label{eq:tempo3}
 p(t(\sigma))=p(t)
\end{equation}
and it follows from the stationary condition of
Theorem~\ref{theo:pmpdfc} (see third equality in
\eqref{eq:dedELfPMP}) that
\begin{equation}
\label{eq:tempo4} \begin{split}
\psi=&-\frac{\partial
\bar{{\mathcal H}}_{f}}{{\partial t_{\sigma}^{'}}} \,
=\frac{\partial}{\partial
t_{\sigma}^{'}}\left[p(t(\sigma))\cdot
\frac{(t_{\sigma}^{'})^{-\alpha}}{\Gamma(n-\alpha)}
\int_{\frac{a}{(t_{\sigma}^{'})^{2}}}^{\sigma}
(\sigma-s)^{n-\alpha-1}\left(\frac{d}{ds}\right)^{n}q(s)ds\right]t_{\sigma}^{'}\\
&  -{\mathcal H}+p\cdot{_a^CD_t^\alpha q}\\
=&-\alpha p(t(\sigma))\cdot
\frac{(t_{\sigma}^{'})^{-\alpha-1}}{\Gamma(n-\alpha)}
\int_{\frac{a}{(t_{\sigma}^{'})^{2}}}^{\sigma}
(\sigma-s)^{n-\alpha-1}\left(\frac{d}{ds}\right)^{n}q(s)ds
-{\mathcal H}+p\cdot{_a^CD_t^\alpha q}\\
=&-\left({\mathcal H}-(1-\alpha)p\cdot{_a^CD_t^\alpha q}\right)\, .
\end{split}
\end{equation}
Substituting \eqref{eq:tempo3} and \eqref{eq:tempo4} into
\eqref{eq:tempo2}, we obtain the fractional conservation law
\eqref{eq:tn-cap}.
\end{proof}

As a corollary, we obtain the analogous to the main result proved
in \cite{CD:FredericoTorres:2007} for fractional problems of the
calculus of variations in the Riemann-Liouville sense.

\begin{definition}[variational invariance for \eqref{eq:pfccv1}]
\label{def:inv-cap2} Functional \eqref{eq:pfccv1} is said to be
invariant under the one-parameter family of infinitesimal
transformations
\begin{equation}
\label{eq:tinfca}
\begin{cases}
\bar{t} = t + \varepsilon\tau(t,q) + o(\varepsilon) \, ,\\
\bar{q}(t) = q(t) + \varepsilon\xi(t,q) + o(\varepsilon) \, ,\\
\end{cases}
\end{equation}
if and only if
\begin{equation*}
\int_{t_{a}}^{t_{b}} L\left(t,q(t),{_a^CD_t^\alpha q(t)}\right) dt =
\int_{\bar{t}(t_a)}^{\bar{t}(t_b)} L\left(\bar{t},\bar{q}(\bar{t}),
{_a^CD_{\bar{t}}^\alpha \bar{q}(\bar{t})},\right) d\bar{t}
\end{equation*}
for any subinterval $[{t_{a}},{t_{b}}] \subseteq [a,b]$.
\end{definition}

\begin{corollary}[Noether's theorem for fractional problems
of the calculus of variations] \label{theo:tnfc2} If functional
\eqref{eq:pfccv1} is invariant under the family of transformations
\eqref{eq:tinfca}, then
\begin{multline}
\label{eq:tncvc}
 C_{f}\left(t,q,{_a^CD_t^\alpha q}\right) =
\partial_{3} L\left(t,q,{_a^CD_t^\alpha q}\right)
 \cdot \xi\\
+ \left[L\left(t,q,{_a^CD_t^\alpha q}\right) - \alpha\partial_{3}
L\left(t,q,{_a^CD_t^\alpha q}\right)\cdot{_a^CD_t^\alpha q}
 \right] \tau
\end{multline}
is a fractional conservation law in the sense of Caputo.
\end{corollary}

\begin{proof}
The fractional conservation law \eqref{eq:tncvc} is obtained
applying Theorem~\ref{theo:tndfc} to functional \eqref{eq:pfccv1}.
\end{proof}

\begin{remark}
If $\alpha=1$ problem \eqref{eq:pfccv1} is reduced to the
classical problem of the calculus of variations,
\begin{equation}
\label{eq:pbcv}
 I[q(\cdot)] = \int_a^b
L\left(t,q(t),\dot{q}(t)\right) \longrightarrow \min \, ,
\end{equation}
and one obtains from Corollary~\ref{theo:tnfc2} the standard
Noether's theorem \cite{Noether:1918}:
\begin{equation}
\label{eq:classNoetCL}
C(t,q,\dot{q})=\partial_{3}
L\left(t,q,\dot{q}\right)
 \cdot \xi(t,q)
+ \left[L\left(t,q,\dot{q}\right) -\partial_{3}
L\left(t,q,\dot{q}\right)\cdot\dot{q}
 \right] \tau(t,q)
\end{equation}
is a conservation law, \textrm{i.e.} \eqref{eq:classNoetCL} is
constant along all the solutions of the Euler-Lagrange equations
\begin{equation}
\label{eq:EL}
\partial_2 L\left(t,q,\dot{q}\right)=\frac{d}{dt}\partial_3 L\left(t,q,\dot{q}\right)
\end{equation}
(these classical equations are obtained from \eqref{eq:E-Lca}
putting $\alpha=1$).
\end{remark}


\section{Illustrative Examples}

In classical mechanics, when problem \eqref{eq:pbcv} does not
depend explicitly on $q$, \textrm{i.e.} $L =
L\left(t,\dot{q}\right)$, it follows from \eqref{eq:dedELfPMP} and
\eqref{eq:EL} that the generalized momentum $p$ is a conservation
law. This is also an immediate consequence of Noether's theorem
\cite{Noether:1918}: from the invariance with respect to
translations on $q$ ($\tau = 0$, $\xi = 1$), it follows from
\eqref{eq:classNoetCL} that $p = \partial_{3} L$ is a conservation
law. Another famous example of application of Noether's theorem in
classical mechanics is given by the conservation of energy: when
the Lagrangian $L$ in \eqref{eq:pbcv} is autonomous, \textrm{i.e.}
$L = L\left(q,\dot{q}\right)$, we have invariance under
time-translations ($\tau = 1$, $\xi = 0$) and it follows from
\eqref{eq:Hnd}, \eqref{eq:dedELfPMP} and \eqref{eq:classNoetCL}
that the Hamiltonian ${\mathcal H}$ (which is interpreted as being the
energy in classical mechanics) is a conservation law. Surprisingly
enough, we show next, as an immediate consequence of our
Theorem~\ref{theo:tndfc}, that for the problem \eqref{Pcap} with a
fractional order of differentiation $\alpha$ ($\alpha \ne 1$), the
following happens:

\begin{itemize}

\item[(i)] similarly to classical mechanics, the generalized
momentum $p$ is a fractional conservation law when $L$ and
$\varphi$ do not depend explicitly on $q$
(Example~\ref{ex:consMom});

\item[(ii)] differently from classical mechanics, the Hamiltonian
${\mathcal H}$ is not a fractional conservation law when $L$ and
$\varphi$ are autonomous (Example~\ref{ex:consEng}).
\end{itemize}
In situation (ii), we obtain from our Theorem~\ref{theo:tndfc} a
new fractional conservation law that involves not only the
Hamiltonian ${\mathcal H}$ but also the fractional order of
differentiation $\alpha$, the generalized momentum $p$, and the
Caputo derivative of the state trajectory $q$ (see
\eqref{eq:eneca} below). This is in agreement with the claim that
the fractional calculus of variations provide a very good
formalism to model nonconservative mechanics
\cite{El-NabulsiMMAS,CD:Riewe:1996}. In the classical case we have
$\alpha = 1$ and the new obtained fractional conservation law
\eqref{eq:eneca} reduces to the expected ``conservation of
energy'' ${\mathcal H}$.

\begin{example}
\label{ex:consMom} Let us consider problem \eqref{Pcap} with
$L(t,q,u)=L(t,u)$, $\varphi(t,q,u)=\varphi(t,u)$. Such a problem
is invariant under translations on the variable $q$, \textrm{i.e.}
condition \eqref{eq:invnd-co} is verified for $\bar{t}=t$,
$\bar{q}(t)=q(t)+\varepsilon$, $\bar{u}(\bar{t})=u(t)$ and
$\bar{p}(\bar{t})=p(t)$: we have $d\bar{t}=dt$ and condition
\eqref{eq:invnd-co} is satisfied since
${_a^CD^\alpha_{\bar{t}}}\bar{q}(\bar{t})={_a^CD^\alpha_{t}}q(t)$:
\begin{equation*}
\begin{split}
_a^CD^\alpha_{\bar{t}}\bar{q}(\bar{t})
&=\frac{1}{\Gamma(n-\alpha)} \int_{\bar{a}}^{\bar{t}}
(\bar{t}-\theta)^{n-\alpha-1}\left(\frac{d}{d\theta}\right)^{n}
\bar{q}(\theta)d\theta
\\
&=\frac{1}{\Gamma(n-\alpha)} \int_{a}^{t}
(t-\theta)^{n-\alpha-1}\left(\frac{d}{d\theta}\right)^{n}
\left(q(t)+\varepsilon\right)d\theta
\\
&={_a^CD^\alpha_{t}}q(t)+{_a^CD^\alpha_{t}}\varepsilon\\
&={_a^CD^\alpha_{t}}q(t)\, .
\end{split}
\end{equation*}
According with \eqref{eq:tinf} one has $\xi=1$ and
$\tau=\varsigma=\varrho=0$. It follows from
Theorem~\ref{theo:tndfc} that $p(t)$
is a fractional conservation law in the sense of Caputo.
\end{example}

\begin{example}
\label{ex:consEng}
We now consider the autonomous
problem \eqref{Pcap}:
$L(t,q,u)=L(q,u)$ and $\varphi(t,q,u)=\varphi(q,u)$.
This problem is invariant under time translation,
\textrm{i.e.} the invariance condition \eqref{eq:invnd-co}
is verified for $\bar{t}=t+\varepsilon$,
$\bar{q}(\bar{t})=q(t)$, $\bar{u}(\bar{t})=u(t)$ and
$\bar{p}(\bar{t})=p(t)$: we have $d\bar{t}=dt$ and
\eqref{eq:invnd-co} follows from the fact that
${_a^CD^\alpha_{\bar{t}}}\bar{q}(\bar{t})={_a^CD^\alpha_{t}}q(t)$:
\begin{equation*}
\begin{split}
_a^CD^\alpha_{\bar{t}}\bar{q}(\bar{t})
&=\frac{1}{\Gamma(n-\alpha)} \int_{\bar{a}}^{\bar{t}}
(\bar{t}-\theta)^{n-\alpha-1}\left(\frac{d}{d\theta}\right)^{n}
\bar{q}(\theta)d\theta
\\
&=\frac{1}{\Gamma(n-\alpha)} \int_{a+\varepsilon}^{t+\varepsilon}
(t+\varepsilon-\theta)^{n-\alpha-1}\left(\frac{d}{d\theta}\right)^{n}
\bar{q}(\theta)d\theta\\
&=\frac{1}{\Gamma(n-\alpha)} \int_{a}^{t}
(t-s)^{n-\alpha-1}\left(\frac{d}{ds}\right)^{n}
\bar{q}(t+\varepsilon)ds\\
&={_a^CD^\alpha_{t}\bar{q}(t+\varepsilon)}={_a^CD^\alpha_{t}\bar{q}(\bar{t})}\\
&={_a^CD^\alpha_{t}{q}(t)}\, .
\end{split}
\end{equation*}
With the notation \eqref{eq:tinf} one has $\tau=1$ and
$\xi=\varsigma=\varrho=0$. We conclude from Theorem~\ref{theo:tndfc} that
\begin{equation}
\label{eq:eneca} {\mathcal H}(t,q,u,p)-(1-\alpha)p\cdot{_a^CD_t^\alpha q}
\end{equation}
is a fractional conservation law in the sense of Caputo. For $\alpha=1$
\eqref{eq:eneca} represents the ``conservation of the total energy'':
\begin{equation*}
{\mathcal H}(t,q(t),u(t),p(t))= constant \, , \quad
t \in [a,b] \, ,
\end{equation*}
for any Pontryagin extremal $(q(\cdot),u(\cdot),p(\cdot))$
of the problem.
\end{example}


{\bf ACKNOWLEDGEMENTS.} This work is part of the first author's PhD project. It was financially supported by the Portuguese Institute for Development (IPAD, G.F.); and the Centre for Research on Optimization and Control (CEOC, D.T.) through the Portuguese Foundation for Science and Technology (FCT) and the European Social Fund.



\end{document}